\documentclass{IEEEtran4PSCC}
\usepackage{cite}
\ifCLASSINFOpdf
   \usepackage[pdftex]{graphicx}
\else
   \usepackage[dvips]{graphicx}
\fi
\usepackage[cmex10]{amsmath}
\usepackage{amssymb}
\usepackage{algorithmic}
\usepackage{tabularx}
\usepackage{color}
\usepackage{soul}

\begin{document}
\title{From zonal to nodal capacity expansion planning: Spatial aggregation impacts on a realistic test-case}
\author{\IEEEauthorblockN{Elizabeth Glista\IEEEauthorrefmark{1},
Bernard Knueven\IEEEauthorrefmark{2},
and Jean-Paul Watson\IEEEauthorrefmark{1}}
}

\maketitle

\begin{abstract}
 Solving power system capacity expansion planning (CEP) problems at realistic spatial resolutions is computationally challenging. Thus, a common practice is to solve CEP over zonal models with low spatial resolution rather than over full-scale nodal power networks.  Due to improvements in solving large-scale stochastic mixed integer programs, these computational limitations are becoming less relevant, and the assumption that zonal models are realistic and useful approximations of nodal CEP is worth revisiting. This work is the first to conduct a systematic computational study on the assumption that spatial aggregation can reasonably be used for ISO-scale CEP. By considering a realistic, large-scale test network based on the state of California with over 8,000 buses, we find that well-designed small spatial aggregations can yield good approximations but that coarser zonal models may result in large distortions of investment decisions, e.g., capacity under-investment of up to 41\% for the lowest resolution model considered. 
\end{abstract}


\thanksto{
\indent\IEEEauthorrefmark{1}Elizabeth Glista and Jean-Paul Watson are with Lawrence Livermore National Laboratory (LLNL). Emails: \{glista1; watson61\}@llnl.gov. This work was performed under the auspices of the U.S. Department of
Energy by Lawrence Livermore National Laboratory under Contract DE-AC52-07NA27344. The authors gratefully acknowledge Gurobi Optimization for providing a free academic license for their software, which was used in this work.
\\
\indent \IEEEauthorrefmark{2}Bernard Knueven is with the National Laboratory of the Rockies. Email: bernard.knueven@nlr.gov. This work was authored in part by the National Laboratory of the Rockies for the U.S. Department of Energy (DOE), operated under Contract No. DE-AC36-08GO28308.
\\
\indent This work was supported by the Department of Energy (DOE)
Office of Electricity’s Advanced Grid Modeling program and the LLNL-LDRD Program under Project No. LDRD25-SI-007.  
The views expressed in the article do not necessarily represent the views of the DOE or the U.S. Government. The U.S. Government retains and the publisher, by accepting the article for publication, acknowledges that the U.S. Government retains a nonexclusive, paid-up, irrevocable, worldwide license to publish or reproduce the published form of this work, or allow others to do so, for U.S. Government purposes. LLNL-CONF-2012972.
}

\section{Introduction}
Due to booming demand for electric power and an aging electric infrastructure, new investments are needed in most existing power systems.  The power system capacity expansion planning (CEP) problem aims to find a system-optimal investment portfolio of new generation, storage, and transmission. CEP is a well-studied problem in the power systems research literature \cite{conejo_investment_2016,gacitua_comprehensive_2018}, and various papers have demonstrated how different modeling choices can improve solution quality.  Some of these include structural choices like co-optimizing generation, storage, and transmission \cite{Go2016} or using stochastic models to capture long-term uncertainty \cite{Scott2021, munoz_2014}.  The choice to include certain features of operational detail has been studied in relation to unit commitment constraints \cite{palmintier_impact_2016,poncelet_unit_2020}, storage modeling \cite{levin_energy_2023,ouwerkerk_impacts_2022}, and power flow representation \cite{ouwerkerk_impacts_2022,neumann_assessments_2022,xu_value_2019}.  Higher temporal resolution, allowing for higher fidelity modeling of intermittent power sources, was studied in \cite{mallapragada_impact_2018,marcy_comparison_2022} and higher spatial resolution in \cite{krishnan_evaluating_2016,frysztacki_strong_2021}.  It is impossible to reconcile all these various studies to definitively say that one type of higher resolution-- whether structural, operational, temporal, or spatial-- is best for every power system.  Nonetheless, some work has been done to systematically compare the various assumptions \cite{xu_value_2019,frew_temporal_2016,Jacobson_2024}.  In \cite{xu_value_2019}, the authors found that capturing long-term structural uncertainty was paramount.  
For comparing the interactions of operational, temporal, and spatial resolution, the recent work \cite{Jacobson_2024} provides a systematic comparison, finding that spatial resolution was the most critical dimension for their test case.  

However, the vast majority of existing CEP studies use either zonal models or unrealistically small nodal models, without considering the underlying large-scale nodal power networks.  For example, in \cite{Jacobson_2024}, the largest model considered is a 26-zone model of the United States, and the authors do not validate their results on more realistic nodal models for power system operation.  Many papers like \cite{xu_value_2019} use small nodal models, on the order of 100 buses, to represent wide-area (interconnection-level) power networks that have thousands to tens of thousands of buses.  

While the choice to approximate large-scale nodal power networks with coarse zonal models has historically made sense due to computational limitations, this is becoming less true.  Recent improvements in decomposition techniques, such as progressive hedging (PH) \cite{mpi-sppy} and Benders (L-shaped) decomposition \cite{li_mixedinteger_2022}, enable solving large-scale, mixed integer stochastic programs in reasonable timeframes for long-term planning models.  Using a parallel implementation of PH on HPC clusters, high-resolution CEP problems with thousands of buses and hundreds of representative days have been solved to low optimality gaps in several hours \cite{zuluaga_parallel_2024,glista_electricity_2024}.

Building on these recent computational developments, we conduct a systematic study on the foundational assumption that zonal CEP models are realistic and useful approximations of the underlying nodal CEP models.  In this work, we propose a heuristic spatial reduction method that reduces the dimensionality of a power network.  We show how our method can be used to generate ``collapsed'' variations of a power system, up to a coarse resolution corresponding to a zonal model.  To compare the quality of the various models, we propose a two-step approach that allows for a mapping of results from a low-resolution model back to a high-resolution model.  Starting with an 8,795-bus synthetic-but-realistic model of California based on \cite{Taylor2024}, we consider several variations of the model with lower spatial resolutions, down to a version with 244 buses.  

\section{Problem Background}
We consider a stochastic, two-stage nodal CEP model based on \cite{Go2016} that co-optimizes generation, transmission, and storage investments in a power system.  The model is formulated as a stochastic mixed integer linear program (MILP) with both investment-stage and operational-stage costs and constraints.  Investment-stage constraints include limits on how much land is available for new generators and how much new capacity can be built.  Operational-stage constraints include those for power balance, intermittent generation availability by hour, daily hydropower availability, storage operation, transmission line thermal limits, and renewable portfolio standard (RPS) goals.  We model stochasticity in generation availability and electric power demand and incorporate this via representative days.  We allow for load shedding and violations of the RPS but discourage these via penalty terms.  For the results in Section \ref{sec:simulations}, we neglect a power flow representation like Btheta or PTDF for computational tractability, instead considering a pipe-and-bubble or ``transportation'' model, but we consider transmission line losses based on the formulation in \cite{zuluaga_parallel_2024} for some of the simulations.  For transmission expansion, we consider reinforcing existing lines in the system, e.g., via reconductoring, so we have a binary decision variable associated with each existing line in the network.  For generation expansion, we consider some types of new generation as integer decisions and some as continuous, based on the scale of per unit capacity, e.g., a new nuclear power plant corresponds to an integer decision and a new solar photovoltaic plant can be approximated by a continuous decision variable.  For a complete description of the modeling details, see \cite{Musselman2024}.  Note that our model formulation is equivalent to the one in \cite{Musselman2024}, with the omission of Btheta constraints for computational tractability as described above.

It should be noted that CEP models (both zonal and nodal) are most commonly found in the research literature and/or used to conduct sensitivity studies, e.g., to evaluate different policy scenarios as the U.S. Congressional Budget Office (CBO) does \cite{cbo_report}.  In contrast, existing ISO and regional planning processes, such as CAISO's transmission planning process \cite{CAISO_TPP_2024_2025}, often use nodal network representations with a defined set of study assumptions, e.g., load growth scenarios as in WECC's report \cite{WECC_Year20_CompoundLoadImpacts_2024}, to conduct operational performance analyses, e.g., via Production Cost Modeling (PCM).  These processes are primarily evaluative and typically do not involve solving an integrated stochastic optimization problem that co-optimizes generation, storage, and transmission investments as we do in CEP.

\section{Power System Network Reduction}\label{sec:reductions}
While much of the CEP literature considers zonal rather than nodal models, the literature is relatively limited on how zonal CEP models can be mapped back to underlying nodal models.  For example, the commonly used CEP model ReEDS uses 134 regions to represent the entire CONUS area, modeling power transfers across roughly 300 interzonal interfaces and using heuristics to estimate intraregional transmission upgrade costs \cite{reeds}. While ReEDS outputs can be used as inputs to a \textit{zonal} PCM, there is limited information about how these zonal outputs can be translated to the underlying nodal power network.  The closest work that addresses the zonal mapping issue is \cite{lumbreras_largescale_2017}, in which the authors present a strategy for constructing zonal models for transmission expansion planning (TEP) and mapping results back to the nodal model.  However, the largest network they consider has 1,025 buses, and they do not directly solve the full TEP problem on this scale, so there is no validation of their network reduction technique.

In the broader power systems literature, there have been a variety of techniques proposed for power system network reduction.  Techniques like the Kron reduction can exactly preserve network injections in the reduced network but result in dense networks that are not suitable for transmission expansion problems \cite{dorfler_kron_2013}.  More relevant to CEP planning are zonal approximation techniques for market modeling, based on electrical distance as in \cite{blumsack_defining_2009}, power transfer distribution factors (PTDFs) as in \cite{lumbreras_largescale_2017,shi_improved_2012}, or locational marginal prices (LMPs) or congestion prices as in \cite{singh_2005,stockl_2025}.  These zonal techniques may rely on pre-determined zones, e.g., in \cite{lumbreras_largescale_2017} the authors propose a measure that combines geographical and electrical distances to define the different zones.  In the next section, we propose a straightforward heuristic method that can be used to define zones from a high-resolution nodal model.  Our heuristic uses geographic distance as the metric for finding nearby buses to merge, as this best preserves geographically localized information about loads and intermittent generation availability that is pertinent to our CEP application and model formulation.  However, electrical distance could also be used and may be more suitable when incorporating more realistic power flow representations in the model.  For the purposes of our study on CEP using a pipe-and-bubble transmission model, with or without losses, this aggregation heuristic captures the overall features of generation availability and load location and is conservative with respect to interzonal flow limitations.  This proposed network reduction approach is just one possible approach and is not necessarily representative of how the entire literature of network reduction techniques would perform on a CEP model.  Future work on CEP using more advanced transmission models may benefit from combining our methodology with other network reduction approaches in the literature, e.g., the approach in \cite{lumbreras_largescale_2017}.

\subsection{Proposed Network Reduction Method}\label{sec:our_method}
Let a full-resolution power system be defined by a set of buses $\mathcal{B}$, transmission lines $\mathcal{L}$, transformers $\mathcal{X}$, generators $\mathcal{G}$, storage elements $\mathcal{S}$, and loads $\mathcal{A}$. Let $o(\ell)$ and $d(\ell)$ represent the ``from'' and ``to'' buses for branches $\ell\in\mathcal{L}\cup\mathcal{X}$. Let $dist(i,j)$ represent the geographic distance (obtained from the Haversine formula) between buses $i,j\in\mathcal{B}$.  We use the notation $k_i$ to denote the degree of bus $i\in\mathcal{B}$, i.e., the number of branches connected to the bus.  The set of radial lines is the set $\mathcal{L}^r\triangleq  \{l\in\mathcal{L}:k_{o(\ell)}=1 \vee k_{d(\ell)}=1\}$. We use the pi-model for transmission lines in the network, in which the series impedance of line $\ell$ is given by $r_\ell+\textbf{j}x_\ell$, where $r_\ell$ is the series resistance and $x_\ell$ is the series reactance.

\begin{figure}[t]
    \centering
    \includegraphics[width=2.8in]{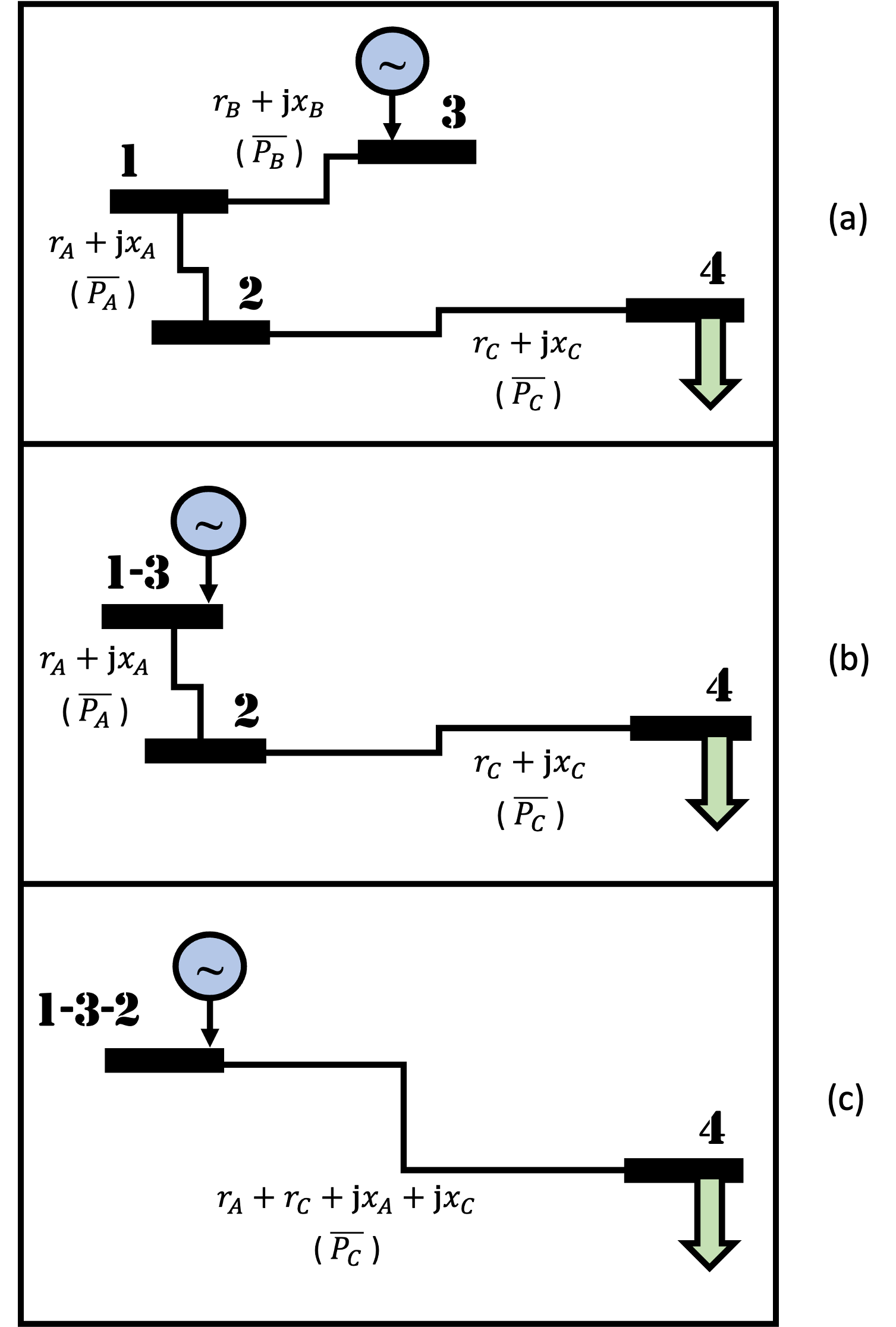}
    \caption{Proposed network reduction method on a simple example.  Transmission line parameters of resistance $r$, reactance $x$, and line rating $\overline{P}$ are shown.  The original 4-bus network is shown in (a).  First, short radial lines are eliminated, as shown in (b).  Then, the remaining bus within the given distance threshold (i.e., bus 2) is merged, as shown in (c).}
    \label{fig:merge_diagram}
\end{figure}

Given some geographic distance threshold $D\geq 0$, we first consider the set of radial transmission lines $\mathcal{L}^r$.  For radial lines with buses within the geographic distance, i.e., the set $\{\ell\in\mathcal{L}^r:dist(o(\ell),d(\ell))\leq D\}$, we remove the bus with degree equal to 1 as well as the radial line.  If there were loads, generators, or storage (either existing or candidates for new investment) at this bus, those are re-located to the remaining bus.  We continue this process until there are no remaining radial lines in the reduced network that have ``from'' and ``to'' buses located within the distance threshold.  Next, we consider the remaining set of transmission lines with buses within the distance threshold.  For each pair of buses that we merge, we define one bus to be the primary bus and one bus to be the secondary bus.  The priority order for assigning the primary bus is to choose (1) the substation bus if only one bus corresponds to a substation; (2) the bus that has already been merged with other buses if the other bus has not been merged; (3) the more connected bus, i.e., the bus with greater degree; (4) the ``from'' bus (arbitrary).  This ordering is a heuristic based on the goal of first preserving the substation buses from the original network in the reduced network and then by prioritizing merged and highly connected buses. The new merged bus adopts the geographic location of the primary bus (including the wind and solar availability at that bus), and adjacent lines are treated as connecting to the primary bus, using electrical equivalence rules to convert transmission line resistances and reactances in series and in parallel.  Transmission line ratings for adjacent lines are preserved. See Figure \ref{fig:merge_diagram} for an example.  Similar to the case with radial lines, we re-locate loads, generators and/or storage to the merged bus (at the location of the primary bus).  We do not collapse branches that correspond to transformers so that bus voltages remain consistent in the reduced network.  Note that because we neglect power flow constraints such as Btheta or PTDF in the current formulation, the line resistances and reactances are only used in the version with transmission line losses.

In this network reduction method, some operational constraints are relaxed. In particular, the internal transmission lines among a set of merged buses are removed from the model, and the set of merged buses is treated as a copperplate.  This assumption is consistent with many existing zonal approaches including that in the widely used ReEDS model \cite{reeds}.  We can see that this heuristic overestimates deliverability and therefore acts as a lower bound on the original problem (with the caveat that generator availability at merged buses should be treated as the minimum for a true lower bound to hold).  Note that the assumption of unlimited power flow between merged buses is not strictly necessary.  A possible enhancement to this approach would be to impose surrogate internal or interface flow limits for merged sets of buses. Similarly, we could preserve information about the original generator availability of the buses (rather than assuming that the merged buses adopt the availability of the primary bus).  These enhanced approaches are left for future work.

For solving CEP problems on these reduced networks and appropriately comparing them to the baseline CEP problem, we consider a two-step approach:
\begin{align*}
   \quad \text{Step (1) }&\text{Solve CEP on the reduced model}\\
   & \text{ to get some solution }\tilde{\textbf{x}}\in\mathbb{R}^{\tilde{n}}\\
   \quad \text{Step (2) }&\text{Map }\tilde{\textbf{x}}\text{ to the baseline CEP model}\\
   &\text{ and re-solve to get final solution }\textbf{x}'\in\mathbb{R}^n
\end{align*}
where we have $\tilde{n}$ investment variables for the reduced model and $n$ investment variables for the baseline model, with $n>\tilde{n}$. Some strategies for Step (2) mappings are discussed in Section \ref{sec:mapping}.  Note that this two-step framework is applicable to either solving the single-day deterministic CEP problems (where Step (1) involves solving a MILP) or the multi-scenario stochastic problems (where Step (1) involves solving a stochastic MILP).  Additionally, the ``two-step'' nature of this model comparison framework should not be confused with the ``two-stage'' nature of our stochastic CEP problem, which captures the different time horizons in investment and operational decisions.

\subsection{Variations on Reduced Networks}\label{sec:variations}
We also consider several variations of the basic reduction method presented above, in order to better preserve information about the original network and thus improve the fidelity of the reduced models.  We consider a variation where only radial lines are collapsed.  We also consider a variation where we tighten the amount of new generation or storage allowed to be built at a candidate site, based on twice the sum of incident existing line ratings (the factor of 2 is due to the fact that the lines are allowed to be reinforced to double their capacity in Step (2)).  This strategy allows for the preservation of some information about transmission capacity in the reduced model, although it neglects information about loads and losses.

\subsection{Projecting Network Reductions onto the Original Network}\label{sec:mapping}
After solving Step (1) of a reduced CEP model, we need to map the new investments in generation, storage, and transmission back to the original full-resolution model.  For any new generation or storage investment in the reduced model, if there is only one possible bus in the original model, either due to land-use constraints or generator availability, then it is trivial to map the new investments back to the original network.  For example, in Figure \ref{fig:merge_diagram}, if solar PV is only allowed to be built at bus 3 and not at buses 1 or 2 due to land use constraints in the original model, then we can directly map solar PV investment at bus 1-3-2 in the reduced network to bus 3 in the original network.

In the case where there are multiple possible buses in the original model that could correspond to a new generation or storage project in the reduced model, e.g., if solar PV investment is allowed at buses 2 or 3 in the above example, we consider two different strategies. One strategy is to distribute the new generation or storage among the possible buses, spread in proportion to the amount of capacity available at each of the buses, and fix these investments (``Map A'').  The second strategy is to fix the total amount of investment (of the given generation or storage type) among the possible buses and allow for the re-optimization of exact locations in Step (2) (``Map B'').  In the above example, if solar investment at bus 1-3-2 in the reduced model is 50 MW and we have equal capacity available for new solar PV at buses 2 and 3, Map A would correspond to fixing solar PV investment at 25 MW for both buses, and Map B would correspond to allowing for the re-optimization of solar PV investments at buses 2 and 3, with the constraint that they sum to 50 MW.  For Map B, the re-optimization of investment locations in Step (2) gives the coarser resolution models an opportunity to adjust to the higher resolution model, thus giving them the best possible opportunity to reach optimal results. 

As an alternative to the geographically-constrained mapping approaches above, we also consider allowing for a full re-optimization of all locations in the original system and only fix the aggregate generation and storage amounts by type (``Map C'').  For example, if the reduced model results in 10 GW of solar PV investment across the entire network, we only fix the total solar PV investment to be 10 GW but allow for a full re-optimization of the exact sizes and locations of the investments.  This final method allows for an analysis of how different high-level choices from the different models can impact system costs and operation.

For transmission line investments in the merged network, we map the components of merged lines back to the original network.  For example, in Figure \ref{fig:merge_diagram}, if we invest in line A-C in the merged network, this investment can be mapped back to investing in lines A and C in the original network.  For lines in the original network that do not appear in the merged network, we allow for a re-optimization of transmission investments.  In Figure \ref{fig:merge_diagram}, we see an example of this at the radial line B, which does not appear in the merged network.

\section{Realistic California Test System}
To test our proposed approach, we consider the California Test System (CATS) developed in \cite{Taylor2024}.  CATS is a synthetic-but-realistic, open-source model of the state of California with 8,870 buses, 10,162 transmission lines, 661 transformers, and 2,092 generators.  
It covers nearly the entire transmission system of the state of California and thus includes the regions operated by CAISO as well as smaller system operators like the Los Angeles Department of Water and Power (LADWP).  CATS was created by combining publicly available data on California's electrical grid with some estimated parameters such as transmission line resistances and reactances.  However, the process to create a connected grid from GIS transmission line data involved some amount of manual tuning \cite{Taylor2024}.  We noticed several disconnected but geographically close buses in CATS that should have been merged per the methodology outlined in \cite{Taylor2024}.  We identified these disconnected buses that were within 0.1 km of each other but electrically far apart (based on graph-distance) and merged the buses.  Additionally, we noticed several nearby buses at substations that appeared to be missing transformers between them, and we added these transformers using a similar approach to that in \cite{Taylor2024}.  We also added missing generators from 2019 Form EIA-860 as described in \cite{Musselman2024}.  After making these changes to our data, we have a version of CATS with 8,795 buses, 10,162 transmission lines, and 777 transformers.  Our system has 2,544 generators corresponding to 80 GW of capacity and 74 storage elements corresponding to 252 MW (728 MWh) of battery storage and 3 GW (33 GWh) of pumped hydroelectric storage.  We consider this to be our baseline high-fidelity version of the California power system.  

Next, we add relevant data for CEP, including costs and constraints for new investments.  In order to model stochastic variability in CEP, we consider 360 representative 24-hour days corresponding to the future weather year 2045.  Obtained from future weather projections of DOE's Energy Exascale Earth System Model (E3SM) California Regionally Refined Model (CARRM) \cite{zhang_leveraging_2023}, which has a 3 km spatial resolution over California, these representative days include data on hourly wind and solar availability, daily hydropower availability, and hourly electric power demand \cite{monteagudo_skillful_2025}.  For the simulations in this work, we consider an RPS target of 60\%.  See \cite{Musselman2024} for more details on data sources for this model.

\begin{table*}[t]
\caption{Characteristics of Test Networks with Varying Spatial Resolution}
\label{tab:network_characteristics}
\centering
\begin{tabular}{|l|l|l|l|l|l|l|l|}
\hline
\textbf{Case name} & \textbf{Distance} &\textbf{Reduce }&\textbf{\# Buses} & \textbf{\# Branches} & \textbf{Minimum } & \textbf{Minimum} &\textbf{Gurobi memory} \\
 & \textbf{for reduction} &\textbf{meshed lines?}&&  &\textbf{resistance ($r$)} &  \textbf{reactance ($x$)}&\textbf{for root LP solve} \\
\hline
CATS& N/A & N/A & 8795 & 10939 & 6.00e-08 & 1.21e-06 &  $\approx$ 3.0 GB\\
KITTENS-r1 & 1 km & No & 7092 & 9217 & 6.00e-08 & 1.21e-06 &  $\approx$ 2.3 GB\\
KITTENS-r3 & 3 km & No &  6676 & 8794 & 6.00e-08 & 1.21e-06 &  $\approx$ 2.2 GB\\
KITTENS-r10 & 10 km & No &  6385 & 8497 & 6.00e-08 & 1.21e-06  &  $\approx$ 2.0 GB\\
KITTENS-r100 & 100 km & No & 6277 & 8385  & 6.00e-08 & 1.21e-06  &  $\approx$ 2.0 GB\\
KITTENS-1 & 1 km  & Yes & 3928 & 5581 & 1.09e-05 & 2.22e-04 &  $\approx$ 1.3 GB\\
KITTENS-3 & 3 km  & Yes & 2233 & 3534 & 4.01e-05 & 4.49e-04 &  $\approx$ 800 MB\\
KITTENS-10 & 10 km & Yes & 907 & 1839 & 4.01e-05 & 2.41e-03 &  $\approx$ 400 MB\\
KITTENS-100 & 100 km& Yes & 244 & 834 & 4.01e-05 & 3.21e-03 &  $\approx$ 120 MB\\
\hline
\end{tabular}
\end{table*}

\section{Reduced Variations of CATS, a.k.a. KITTENS}
From our baseline high-fidelity version of CATS, we consider several distances for network reduction--1 km, 3 km, 10 km, and 100 km--using the approach described in Section \ref{sec:reductions}.  We use the shorthand KITTENS to refer to these reduced variations, where, with a slight abuse of spelling, KITTENS stands for the \underline{\textbf{K}}ollapsed Cal\underline{\textbf{I}}fornia
\underline{\textbf{T}}est Systems for \underline{\textbf{T}}ransmission \underline{\textbf{E}}xpansion with Co-optimized Ge\underline{\textbf{N}}eration and \underline{\textbf{S}}torage Expansion.

We use the notation KITTENS-1 to refer to the reduced version of CATS based on a 1 km nodal proximity distance, including meshed lines in the reduction.  We use the notation KITTENS-r1 to refer to the version where we only reduce radial lines (as described in Section \ref{sec:variations}) based on a 1 km nodal proximity distance. Similarly, we consider 3 km, 10 km, and 100 km distances.  Thus, KITTENS-r1 is the largest network of the reduced models (and most similar to CATS), and KITTENS-100 is the smallest.  Note that the 1 and 3 km KITTENS models are within the spatial resolution of the E3SM CARRM weather data and thus are within the tolerance for generator availability. The reduced KITTENS models naturally have lower memory and solve time requirements for CEP compared to the full CATS model, and the networks with reduced meshed lines exhibit better numerical stability.  While the original network has line resistance and reactance values as low as \texttt{6e-8} and \texttt{1e-6} p.u. respectively, KITTENS-1 has line resistance and reactance values on the order of \texttt{1e-5} and \texttt{1e-4} p.u.  A summary of the different test networks and their characteristics is given in Table \ref{tab:network_characteristics}.

It is worth highlighting that even the lowest resolution model we consider, KITTENS-100, which has 244 buses for the state of California, is much higher resolution than most models in the existing literature.  For example, in \cite{krishnan_evaluating_2016} the authors use only 4 buses in California while modeling a national model with 134 buses.  In \cite{xu_value_2019}, the case study uses a 300-bus system for the entire Western Electricity Coordinating Council (WECC) region, with only a subset of these buses located in California.

\section{Simulations}\label{sec:simulations}
We now conduct a systematic study of how varying the spatial resolution of our CEP model impacts the solve time, the quality of the solution found, and the recommended investments from the model.  While our main goal is to solve a stochastic CEP problem, it is illustrative to first consider how spatial aggregation impacts single-day deterministic CEP problems.  For this, we solve the 360 single-day CEP scenarios in parallel and compare how solutions from the coarser spatial resolutions map back to the full-resolution network, i.e., the original network without spatial aggregation. Then, we solve the stochastic CEP problem with the full year of data via the parallel decomposition method Progressive Hedging (PH) implemented via \texttt{mpi-sppy} \cite{mpi-sppy}.  We examine how these stochastic results vary with the spatial resolution of the model.  For our given test system, the CATS network corresponds to the benchmark CEP model, and the various KITTENS networks (defined in Table \ref{tab:network_characteristics}) correspond to the reduced models.

To compare solution quality, we define a loss metric called the Reduced Model Mapping Loss (RMML) as the normalized difference between the best objective value (incumbent) of the baseline high-resolution model (i.e., CATS) and the best objective value of the lower resolution model (i.e., KITTENS) as mapped back to the high-resolution model (where choices for this mapping are described in Section \ref{sec:mapping}):
\begin{equation}
    RMML \triangleq \frac{f(\textbf{x}')-f(\textbf{x}^*)}{f(\textbf{x}^*)}*100 \qquad(\%)\label{eqn:error}
\end{equation}
where for some variation in model resolution, let $f(\textbf{x}):\mathbb{R}^n\rightarrow \mathbb{R}$ be the evaluation of an investment portfolio $\textbf{x}\in\mathbb{R}^n$ on the high-resolution CEP problem.  Let $\textbf{x}'\in\mathbb{R}^n$ be the investment solution on the low-resolution CEP problem as mapped back to the high-resolution problem, and let $\textbf{x}^*\in\mathbb{R}^n$ be the investment solution on the high-resolution CEP problem, where $\textbf{x}^*\in\text{argmin}f(\textbf{x})$.  This loss metric can be thought of as a type of optimality gap that quantifies how bad our approximations are at solving the underlying CEP problem on the true system.  We will interchangeably call this our ``error'' metric associated with the model reductions as to avoid confusion with the word ``loss,'' which is also used for transmission losses and lost load.

For the stochastic simulations in Section \ref{sec:stochastic_sims}, this loss metric becomes:
\begin{equation}
    RMML \triangleq \frac{\mathbb{E}\left[f(\textbf{x}')\right]-\mathbb{E}\left[f(\textbf{x}^*)\right]}{\mathbb{E}\left[f(\textbf{x}^*)\right]}*100 \qquad(\%)\label{eqn:error_stochastic}
\end{equation}
where the expected value is taken as the probabilistically weighted objective value over all the scenarios, i.e., $\mathbb{E}\left[f(\textbf{x})\right]\triangleq \sum_{\omega\in\Omega}p_\omega f_\omega(\textbf{x})$ for scenarios $\omega\in\Omega$ with probabilities given by $p_\omega$ and scenario-specific objective functions given by $f_\omega(\textbf{x}):\mathbb{R}^n\rightarrow\mathbb{R}$.  Note that this stochastic loss metric is similar to existing metrics in the literature, including the Value of Model Enhancement (VOME) introduced in \cite{xu_value_2019} and used to compare the tradeoffs between including different model features on the solution quality.

To run the simulations, we use the HPC cluster \texttt{dane} at Lawrence Livermore National Laboratory (LLNL).  Each compute node has 112 2.0 GHz cores and 256 GB of RAM, and we fit 18 to 20 problems/node for solving Step (1) problems and 30 problems/node for solving Step (2) problems.

\subsection{Single-day Deterministic CEP Scenarios}\label{sec:single_day_sims}
\begin{figure}
    \centering
    \includegraphics[width=3.4in]{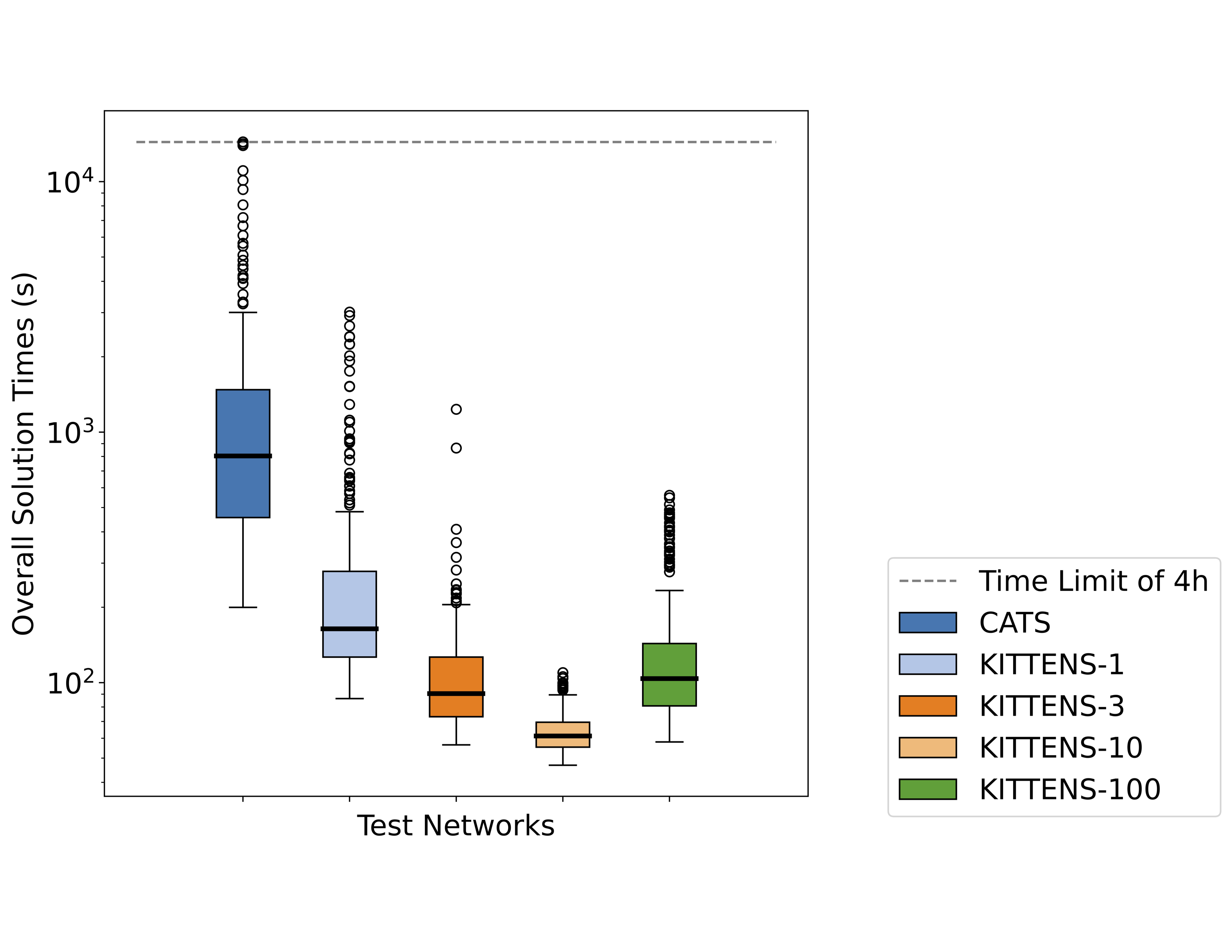}
    \includegraphics[width=3.4in]{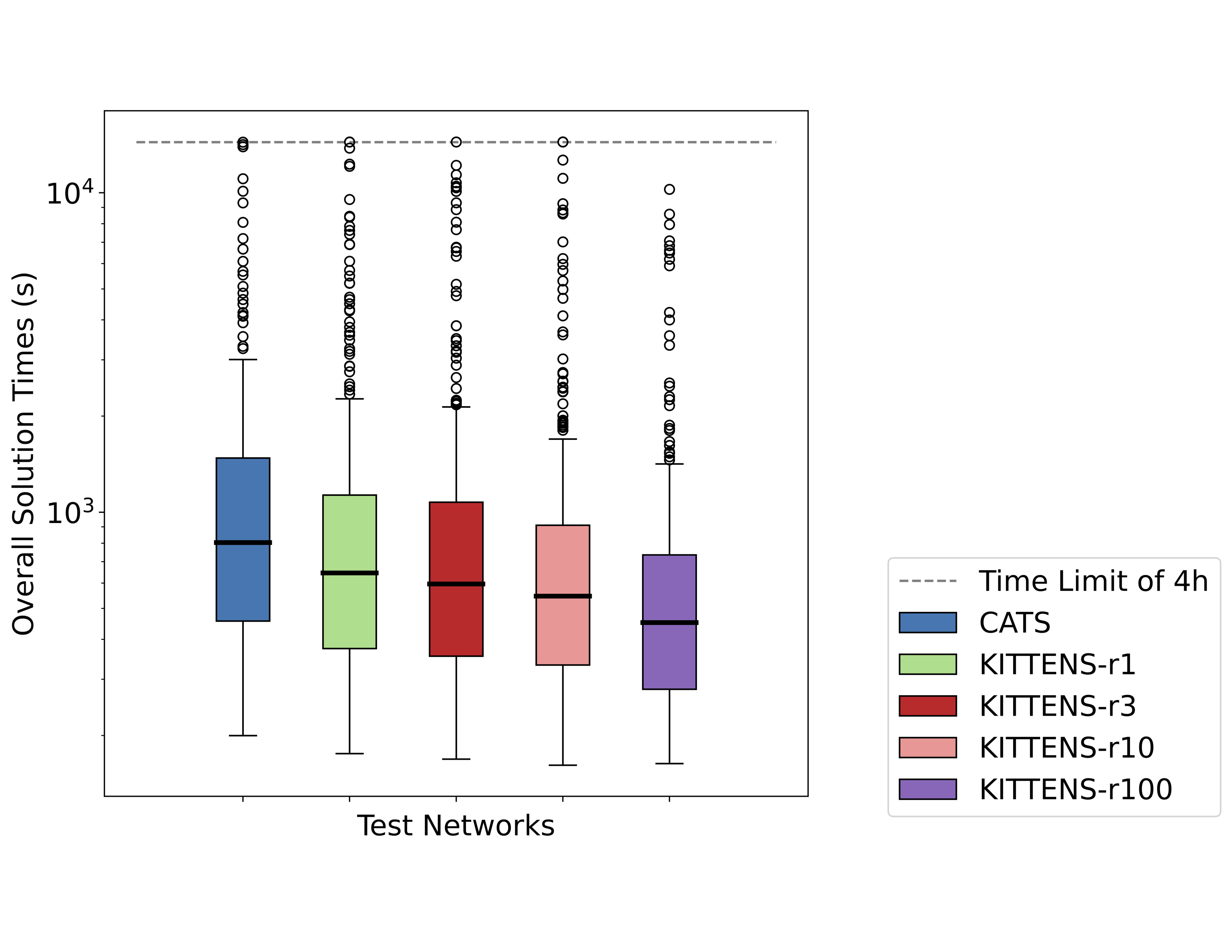}
    \caption{Overall solution time to solve the CEP problem on CATS, with a two-step approach for reduced networks.  Each boxplot corresponds to 360 scenarios of single-day CEP problems.  The top figure includes the KITTENS variations with meshed lines collapsed, and the bottom figure includes the KITTENS variations with only radial lines collapsed.}
    \label{fig:solvetimes_collapsed}
\end{figure}

We solve the 360 single-day CEP scenarios in parallel on \texttt{dane}, using Gurobi to solve the individual sub-problems.  For Step (1), where CEP is solved on the CATS and KITTENS networks, all problems had feasible solutions within the 4-hour time limit, and most but not all solutions were within the target MIP gap of 1\%.  For Step (2), where the solutions from the KITTENS problems are mapped back to CATS and re-solved via Map B (as described in Section \ref{sec:mapping}), all problems solved well before the 1-hour time limit to a MIP gap below the target of 0.1\%.  See Table \ref{tab:first_step_results} for a summary of the Step (1) and (2) results.  The cumulative solve times, combining those in Steps (1) and (2), are visualized in Figure \ref{fig:solvetimes_collapsed}. While solve times decrease quadratically for Step (1) as network size decreases, solve times increase with network size for Step (2) since the mapping has more degrees of freedom when there is less information contained in the Step (1) results.  Thus, we see cumulative solve times decrease until we reach the smallest case, KITTENS-100, at which point they start to increase again. The same simulations with the addition of transmission losses are also included in Table \ref{tab:first_step_results}.  Naturally, the solution times are higher for these simulations compared to the lossless version, but the overall trends between models of varying spatial resolution are the same.  We also consider the ``tightened'' version of the reduced networks, where the amount of allowable generation and storage investment at each bus is reduced based on line capacity in the original network (as described in Section \ref{sec:variations}).  We saw that the tightened versions of the ``radial-only'' KITTENS models (e.g., KITTENS-r1) were noticeably faster to solve than the non-tightened versions; however, this was not true for the ``meshed'' KITTENS models (e.g., KITTENS-1).  On average, solving the two-step CEP model on the tightened KITTENS-r1 was 61\% faster than solving CEP on CATS, and for KITTENS-1 it was 80\% faster.  For some of the more reduced networks, we found that solving CEP on the tightened KITTENS-r100 was 73\% faster than solving it on CATS, and for KITTENS-10 it was 96\% faster.

\begin{figure}
    \centering
    \includegraphics[width=3.4in]{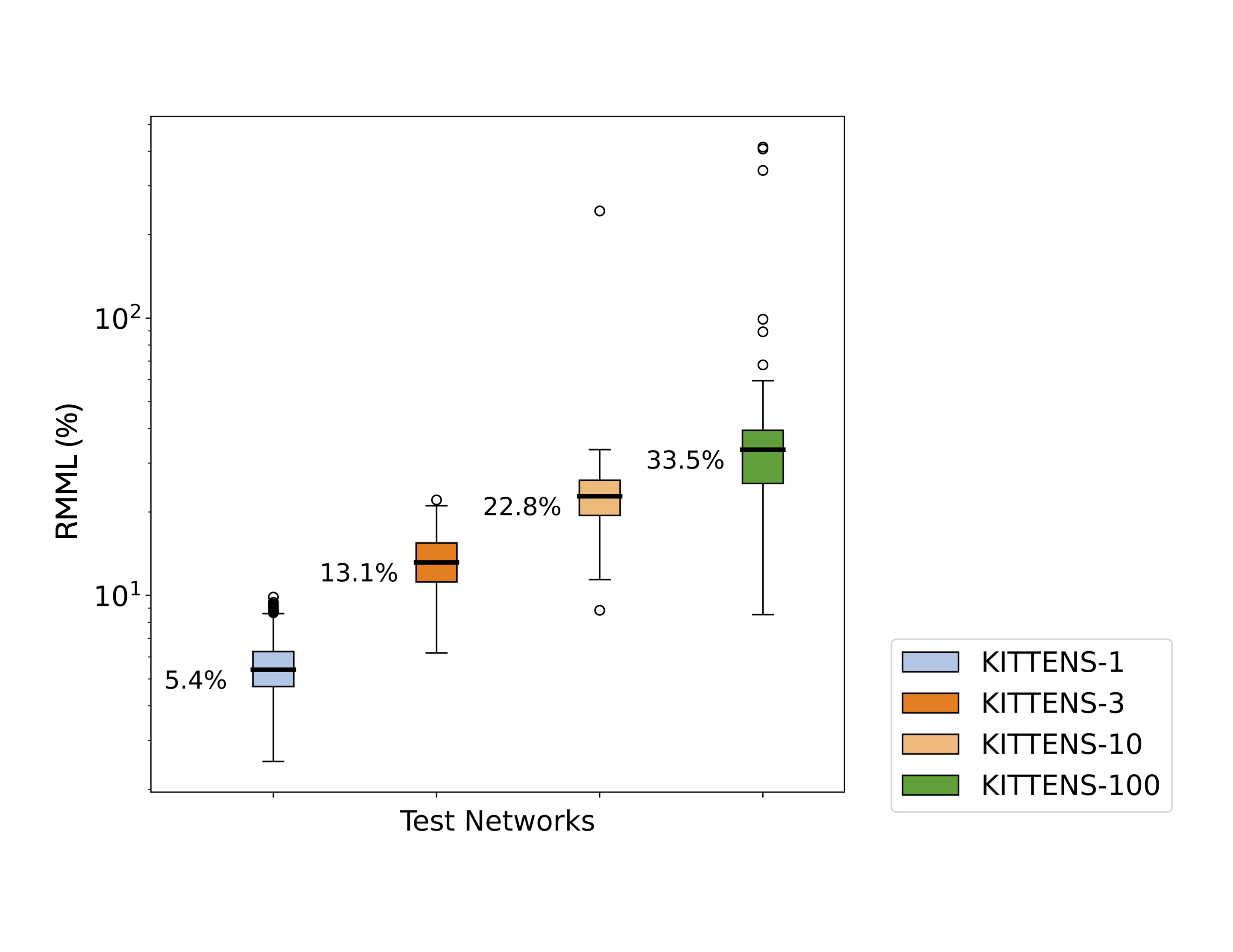}
        \includegraphics[width=3.4in]{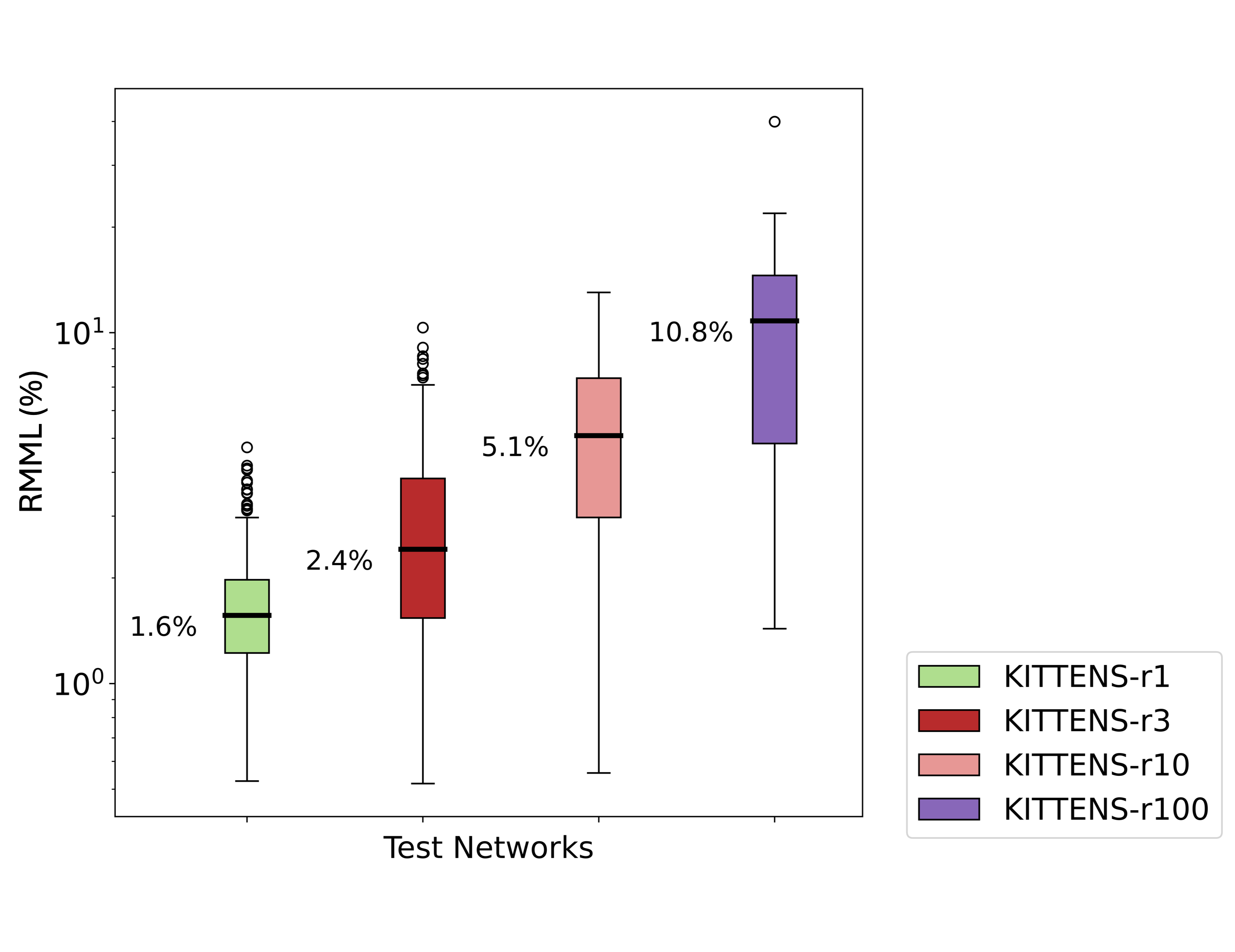}
    \caption{Quality of solutions to reduced KITTENS models when mapped back to CATS, as measured by the error metric RMML given in (\ref{eqn:error}).  Each boxplot corresponds to 360 scenarios of single-day CEP problems.  We see that the majority of reduced models result in median errors greater than 5\%.}
    \label{fig:RMML}
\end{figure}

\begin{figure}
    \centering
    \includegraphics[width=3.4in]{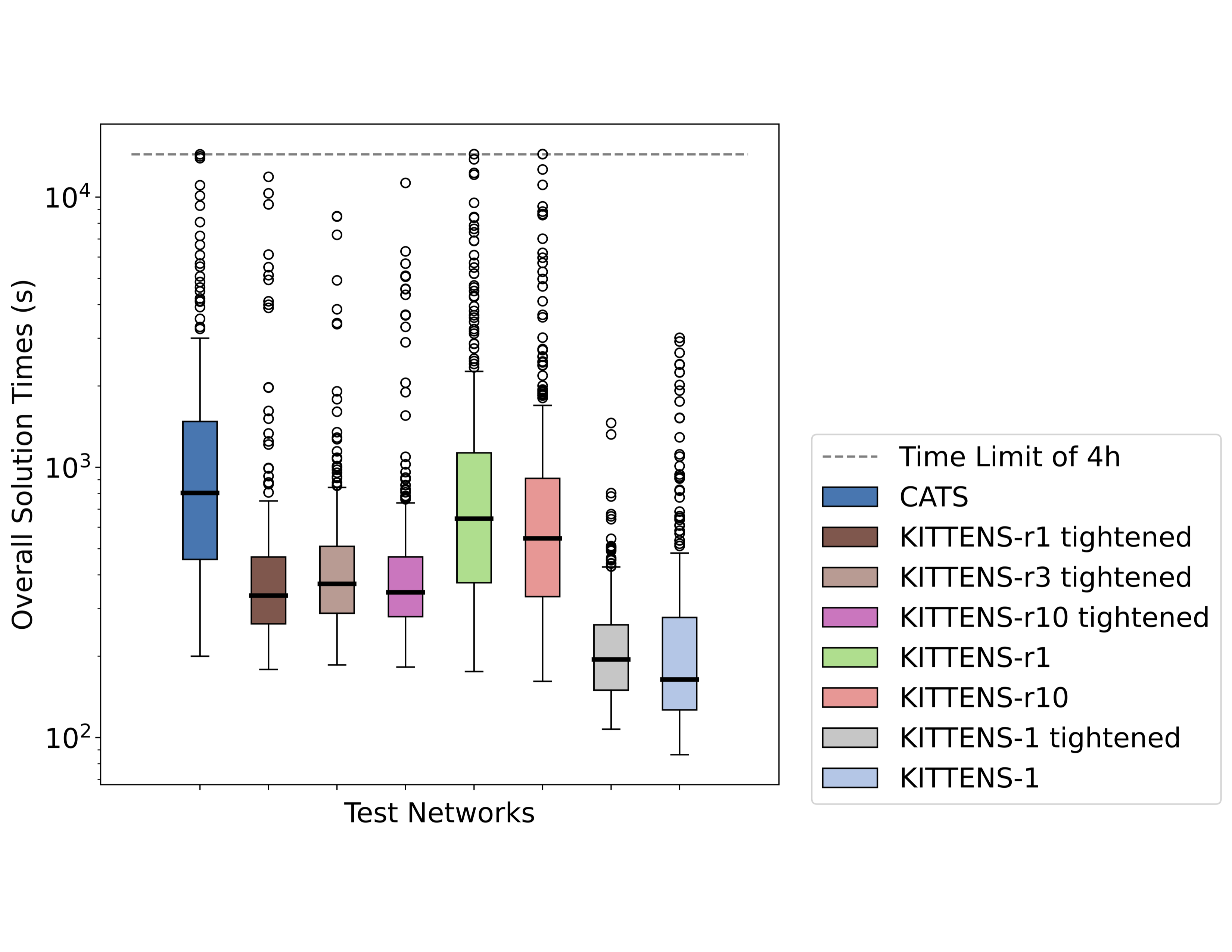}
        \includegraphics[width=3.4in]{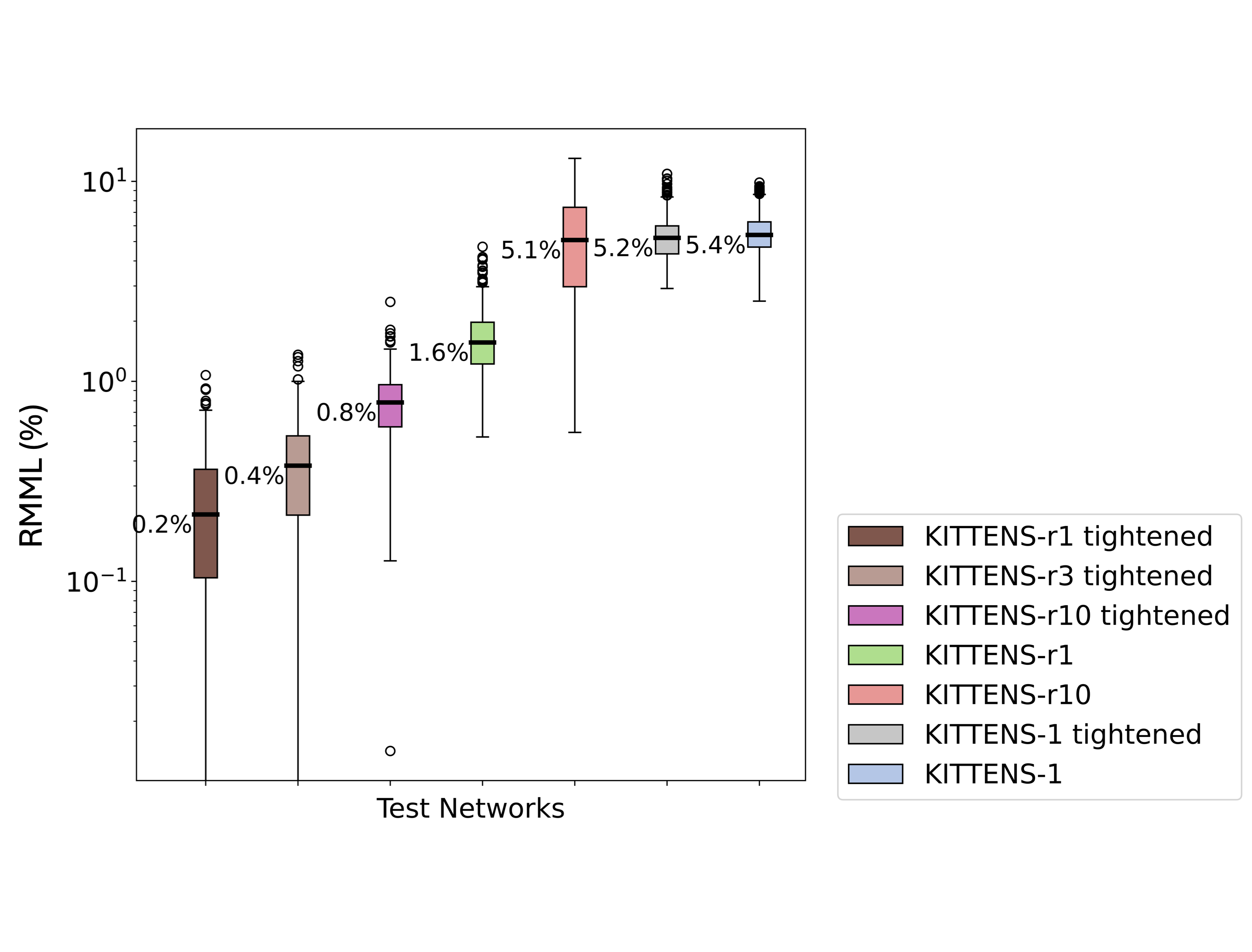}
    \caption{Overall solution times (top figure) and error metric RMML (bottom figure) for a selection of the best performing models, as determined by the median error over all 360 single-day CEP scenarios.  We can see that several of these reduced models obtain median error values below the 1\% optimality tolerance.}
    \label{fig:best_models}
\end{figure}

Comparing the solution quality of the reduced models mapped back to the high-resolution model with the error metric RMML defined in (\ref{eqn:error}), we see that the tightened versions of KITTENS-r1, -r3, and -r10 have average error under 1\%, or within the MIP tolerance (see Table \ref{tab:first_step_results}).  The tightened versions of the radial-only KITTENS models outperform the non-tightened versions, of which only KITTENS-r1 manages to have error under 5\% for all problems.  For the faster meshed KITTENS models, we see that tightening does reduce the error in the majority of problems, although all models still have error over 5\% on average.  The error values for some of these runs are visualized in Figure \ref{fig:RMML} and for the best performing models in Figure \ref{fig:best_models}.  In the case with losses, we see similar values for average error as in the case without losses but much higher values for the maximum error observed over all 360 scenarios.  This implies that low resolution models may be more inadequate in lossy models with congestion compared to lossless or congestion-free scenarios.  Future work will examine how Btheta power flow impacts this analysis, which we hypothesize may yield a larger divergence between the results of models with different spatial resolutions.  

We also re-run the Step (2) problems fixing only the aggregate quantity of new investment of each generation and storage type and allowing for a full re-optimization of the investment locations, corresponding to Map C as described in Section \ref{sec:mapping}.  The Step (2) results of these simulations are shown in Table \ref{tab:second_step_results2}.  Since these simulations allow for more re-optimization than those in Table \ref{tab:first_step_results}, it is natural that all cases have lower average and maximum error values. Similar to the previous results, we see that the tightened KITTENS-r1, -r3 and -r10 models have average error below the 1\% MIP tolerance, with maximum error below 1\% for -r1 and -r3. However, we see that for the more reduced models, e.g., KITTENS-10 and KITTENS-100, we still have high error in many scenarios, with maximum error over 100\% for all but the tightened KITTENS-10.

In Table \ref{tab:buildout_comparison}, we compare how the amount of investment in different technologies changes depending on the spatial representation of the CEP model.  Note that in this table we lump together solar PV and solar thermal (all new investments are in solar PV) and lump together all natural gas (NG) types (CC, CT, ST, ICE).  The category of ``All Other'' contains types found in the existing system that are not considered for investment: conventional hydroelectric, conventional steam coal, petroleum coke, petroleum liquids, imports, and other.  We saw that on average (across all 360 single-day scenarios), for all reduced KITTENS models, there was under-investment in solar PV and over-investment in geothermal, compared to the baseline CATS model.  While smaller in magnitude, we also saw under-investment in storage, and for most KITTENS models, under-investment in onshore wind, landfill gas, and wood biomass. All reduced models resulted in under-investment in new capacity, with the average level of overall under-investment being 0.9\% for KITTENS-r1 ($<$ 0.1\% for tightened version) and 3.4\% for KITTENS-1 (1.3\% for tightened version), up to 41\% for KITTENS-100 (30\% for tightened version).  This under-investment in new generation capacity is likely due to an overestimate of deliverability in the reduced networks.  

\begin{table*}[t]
\caption{Solving Single-Day, Deterministic CEP on Test Networks of Varying Spatial Resolution, via Two-Step Approach}
\label{tab:first_step_results}
\centering
\begin{tabular}{|l|l|l|l|l|l|l|l|l|l|l|}
\hline
&\multicolumn{5}{c|}{\textit{\textbf{Step (1) Problems}}} & \multicolumn{4}{c|}{\textit{\textbf{Step (2) Problems with Map B}}}\\
\hline
 &\textbf{Median} & \textbf{Average} & \textbf{\# of problems} & \textbf{Average} & \textbf{Max.}  & \textbf{Median} & \textbf{Average} & \textbf{Average} & \textbf{Max. error}\\
& \textbf{solve} & \textbf{solve} & \textbf{that hit 4h} & \textbf{MIP} & \textbf{MIP}  & \textbf{solve} & \textbf{solve} & \textbf{error (\%)} & \textbf{(\%) defined} \\
\textbf{Case name}& \textbf{time (s)} & \textbf{time (s)} & \textbf{time limit} &   \textbf{gap (\%)} &  \textbf{gap (\%)} & \textbf{time (s)} & \textbf{time (s)} & \textbf{defined by (\ref{eqn:error})} & \textbf{by (\ref{eqn:error})} \\
\hline
CATS   &   804   &   1436   &   4   &   0.93   &   1.07 & N/A & N/A & N/A & N/A\\
KITTENS-r1   &   631   &   1391   &   6   &   0.93   &   1.00  & 15 & 15 & 1.67 & 4.71 \\
KITTENS-r3   &   580   &   1195   &   2   &   0.93   &   1.08   &16 & 16 & 2.88 & 10.34\\
KITTENS-r10   &   530   &   1090   &   4   &   0.93   &   1.00  & 17 & 17 & 5.40 & 13.03\\
KITTENS-r100   &   435   &   739   &   0   &   0.92   &   1.00  & 18 & 18 & 10.25 & 39.94 \\
KITTENS-1   &   144   &   276   &   0   &   0.92   &   1.00  & 18 & 18 & 5.61 & 9.86\\
KITTENS-3   &   65   &   87   &   0   &   0.91   &   1.00 & 23 & 23 & 13.36 & 22.10 \\
KITTENS-10   &   21   &   23   &   0   &   0.78   &   1.00  & 40 & 41 & 23.31 & 243.66\\
KITTENS-100   &   5   &   5   &   0   &   0.25   &   0.53  & 99 & 133 & 35.62 & 413.86\\
\hline\hline
KITTENS-r1 tightened   &   319   &   548   &   0   &   0.91   &   1.00 & 15 & 15 & \textbf{0.24} & 1.07 \\
KITTENS-r3 tightened   &   353   &   521   &   0   &   0.91   &   1.00 & 16 & 16 & \textbf{0.39} & 1.36  \\
KITTENS-r10 tightened   &   328   &   534   &   0   &   0.91   &   1.00  & 18 & 18 & \textbf{0.79} & 2.50\\
KITTENS-r100 tightened   &   295   &   368   &   0   &   0.92   &   1.00 & 18 & 18 & 4.93 & 22.32  \\
KITTENS-1 tightened   &   180   &   218   &   0   &   0.86   &   1.00 & 13 & 13 & 5.34 & 10.91 \\
KITTENS-3 tightened   &   113   &   130   &   0   &   0.86   &   1.00 & 21 & 22 & 8.73 & 16.42 \\
KITTENS-10 tightened   &   67   &   79   &   0   &   0.73   &   1.00  & 40 & 42 & 15.67 & 77.92\\
KITTENS-100 tightened   &   37   &   40   &   0   &   0.37   &   0.59  & 116 & 166 & 31.29 & 357.05\\\hline\hline
CATS w/losses   &   1104   &   2681   &   12   &   0.70   &   1.74  & N/A & N/A & N/A & N/A\\
KITTENS-r1 w/losses   &   919   &   2243   &   7   &   0.67   &   1.30  & 34 & 36 & 1.30 & 14.26\\
KITTENS-r3 w/losses   &   921   &   2012   &   7   &   0.65   &   1.66  & 25 & 27 & 3.23 & 20.49\\
KITTENS-r10 w/losses   &   833   &   1804   &   4   &   0.65   &   1.16  & 25 & 27 & 5.50 & 32.73\\
KITTENS-r100 w/losses   &   524   &   1531   &   3   &   0.72   &   1.25 & 25 & 26 & 13.24 & 313.53 \\
KITTENS-1 w/losses   &   426   &   923   &   1   &   0.66   &   1.53  & 36 & 39 & 4.80 & 23.96\\
KITTENS-3 w/losses   &   183   &   427   &   0   &   0.72   &   1.00 & 36 & 38 & 12.63 & 94.00 \\
KITTENS-10 w/losses   &   30   &   36   &   0   &   0.79   &   1.00 & 54 & 78 & 28.27 & 674.27 \\
KITTENS-100 w/losses   &   5   &   5   &   0   &   0.27   &   0.66 & 172 & 337 & 48.48 & 1325.99 \\
\hline
\end{tabular}
\end{table*}

\begin{table*}[t]
\caption{Step (2) Results Fixing Only Aggregate Generation and Storage Investment Quantities}
\label{tab:second_step_results2}
\centering
\begin{tabular}{|l|l|l|l|l|l|l|}
\hline
&\multicolumn{4}{c|}{\textit{\textbf{Step (2) Problems with Map C}}}\\
\hline
 &\textbf{Median} & \textbf{Average} & \textbf{Average error (\%)} & \textbf{Max. error (\%)} \\
 \textbf{Case name}& \textbf{solve time (s)} & \textbf{solve time (s)} &   \textbf{defined by (\ref{eqn:error})} & \textbf{defined by (\ref{eqn:error})} \\
 \hline
KITTENS-r1   &   45   &   70   &   \textbf{0.44}   &   1.33  \\
KITTENS-r3   &   49   &   73   &   \textbf{0.69}   &   2.01  \\
KITTENS-r10   &   30   &   60   &   1.24   &   3.63  \\
KITTENS-r100   &   18   &   28   &   5.90   &   16.16  \\
KITTENS-1   &   22   &   46   &   4.46   &   8.73  \\
KITTENS-3   &   30   &   47   &   11.17   &   18.74  \\
KITTENS-10   &   49   &   62   &   17.59   &   189.08  \\
KITTENS-100   &   126   &   202   &   27.55   &   377.09  \\
\hline\hline
KITTENS-r1 tightened   &   25   &   70   &   \textbf{0.08}   &   \textbf{0.89}  \\
KITTENS-r3 tightened   &   40   &   70   &   \textbf{0.13}   &   \textbf{0.91}  \\
KITTENS-r10 tightened   &   39   &   75   &   \textbf{0.22}   &   1.08 \\
KITTENS-r100 tightened   &   19   &   35   &   3.99   &   12.88  \\
KITTENS-1 tightened   &   21   &   40   &   4.68   &   8.56  \\
KITTENS-3 tightened   &   37   &   60   &   7.96   &   16.40  \\
KITTENS-10 tightened   &   58   &   97   &   13.08   &   48.73 \\
KITTENS-100 tightened   &   142   &   247   &   21.92   &   318.81   \\
\hline
\end{tabular}
\end{table*}

\begin{table*}[t]
    \centering
        \caption{Differences in Generation and Storage Investment Choices (in MW) between CEP Models of Varying Spatial Resolution,\\ Averaged Across 360 Single-Day CEP Scenarios }
    \begin{tabular}{|l|l|l|l|l|l|l|l|l|l|l|l|l|l|l|l|}
    \hline
\textbf{Case name} & \textbf{Wood} & \textbf{Landfill} & \textbf{Onshore} & \textbf{Solar}  & \textbf{NG} & \textbf{Nuclear} & \textbf{Geo-} & \textbf{Battery} & \textbf{Hydro}  &\textbf{All}& \textbf{Total} \\
& \textbf{Biomass} & \textbf{Gas} & \textbf{Wind} & \textbf{(all)}  & \textbf{(all) }& & \textbf{thermal} &  &   \textbf{Storage} &\textbf{Other} &\\
\hline\hline
Existing CATS & 648 & 533 & 6163 & 12614 & 39711 & 2323 & 2731 & 253 & 3312& 15274 & 83562\\
\hline\hline
CATS new investments & 345 & 105 & 1136 & 20601  & 1 & 0  & 3250 & 16 & 99& 0 &25553\\
\hline\hline
KITTENS-r1, $\Delta$ invest & -25 & -3 & -28 & -288  & 0 & 0 & 111 & 3& -6& 0 &-236\\
KITTENS-r3, $\Delta$ invest & -36 & -4 & -36 & -470  & 1 & 0 & 110 & -1 & -4& 0 &-440\\
KITTENS-r10, $\Delta$ invest & -68 & -9 & -51 & -670   & 0 & 0 & 109 & -4 & -7& 0 & -700\\
KITTENS-r100, $\Delta$ invest & -255 & -46 & -257 & -3997  & 0 & 0 & 243 & -11 & -55& 0& -4380\\
KITTENS-1, $\Delta$ invest & -88 & -10 & -76 & -1004  & -1 & 0  & 322 & -1  & -15& 0 &-873\\
KITTENS-3, $\Delta$ invest & -154 & -20 & -128 & -1782 & 0 & 0  & 322 & -8 & -38&0 & -1807\\
KITTENS-10, $\Delta$ invest & -273 & -49 & -234 & -2648  & -1 & 0  & 318 & -14 & -82& 0 &-2982\\
KITTENS-100, $\Delta$ invest & -331 & -98 & -654 & -9904  & -1 & 0  & 538 & -16 & -99& 0 & -10563\\
\hline\hline
KITTENS-r1 tight, $\Delta$ invest & 4 & 0 & 2 & -68  & -1  & 0 & 46 & -3 & 0 & 0& -19\\
KITTENS-r3 tight, $\Delta$ invest & 6 & 0 & 2 & -225  &0 & 0 & 49 & -4 & 1& 0 &-172\\
KITTENS-r10 tight, $\Delta$ invest & 4 & 1 & -16 & -346   & 1 & 0 & 51 & -5 & -3&0 & -314\\
KITTENS-r100 tight, $\Delta$ invest & -203 & -32 & -212 & -3313   & -1 & 0 & 216 & -11 & -54& 0 &-3611\\
KITTENS-1 tight, $\Delta$ invest & -33 & -4 & -33 & -490  &  0 & 0 & 231 & -7 & -3& 0 &-340\\
KITTENS-3 tight, $\Delta$ invest & -25 & -3 & -18 & -499 & 0 & 0 & 142 & -8 & -3& 0 &-415\\
KITTENS-10 tight, $\Delta$ invest & -175 & -23 & -113 & -768  &  -1  & 0 & 130 & -11 & -75& 0 &-1035\\
KITTENS-100 tight, $\Delta$ invest & -319 & -85 & -572 & -6783  & -1  & 0 & 315 & -15 & -96&0& -7555\\
\hline
    \end{tabular}
    \vspace{1em}
    \label{tab:buildout_comparison}
\end{table*}

\begin{table*}[t]
    \centering
        \caption{Differences in Generation and Storage Investment Choices (in MW) between CEP Models of Varying Spatial Resolution,\\ For 360-day Stochastic CEP Models }
    \begin{tabular}{|l|l|l|l|l|l|l|l|l|l|l|l|l|l|l|l|}
    \hline
\textbf{Case name} & \textbf{Wood} & \textbf{Landfill} & \textbf{Onshore} & \textbf{Solar}  & \textbf{NG} & \textbf{Nuclear} & \textbf{Geo-} & \textbf{Battery} & \textbf{Hydro}  &\textbf{All}& \textbf{Total} \\
& \textbf{Biomass} & \textbf{Gas} & \textbf{Wind} & \textbf{(all)}  & \textbf{(all) }& & \textbf{thermal} &  &   \textbf{Storage} &\textbf{Other} &\\
\hline\hline
Existing CATS & 648 & 533 & 6163 & 12614 & 39711 & 2323 & 2731 & 253 & 3312& 15274 & 83562\\
\hline\hline
CATS new investments & 711 & 237 & 825 & 30994 & 111 & 0 & 3246 & 91 & 776 & 0 & 36990\\
\hline\hline
KITTENS-r3 tight, $\Delta$ invest & -449 & 1 & -196 &  -2351 & -111 &0 & 0 & 674 & -529 & 0 & -2961 \\
KITTENS-1 tight, $\Delta$ invest & -406 & 1 & -299 & -1622 & 0 & 0 & 160 & -68 & -41 & 0 & -2275 \\
KITTENS-10 tight, $\Delta$ invest & -699 & 1 & -777 & -7070 & -111 & 0 & 0 & -88 & -776 &0 & -9520 \\
KITTENS-100 tight, $\Delta$ invest & -711 & 1 & -825 & -11139 & -111 & 0 & 192 & -91 & -776 & 0 & -13459\\

\hline
    \end{tabular}
    \vspace{1em}
    \label{tab:buildout_comparison_stochastic}
\end{table*}

\subsection{Stochastic CEP on a Full Year of Data}\label{sec:stochastic_sims}
Next, we solve stochastic CEP problems across all 360 representative days.  Because these problems are more computationally intensive to solve, we consider only 4 versions of KITTENS (-r3, -1, -10, and -100, all tightened) in addition to CATS.  We set a time limit of 20 hours for the PH algorithm and obtained optimality gaps around or below 3\% for all KITTENS models and a gap of about 5\% for the CATS model. 
Because re-optimizing generation and storage investments in Step (2) is less straightforward if you have to enforce consensus among the different scenarios, we considered Map A from Section \ref{sec:mapping}, which directly maps generation and storage investments obtained from PH on KITTENS to the CATS model.  For transmission investments, we first considered the approach of mapping components of merged lines back to the original network and then investing in all lines that do not appear in the merged network.  This resulted in good performance for KITTENS-r3, yielding an error of around 1\%, but poor performance for KITTENS-1, -10, and -100.  We also saw counterintuitive behavior where the coarser KITTENS-100 model performed better than the KITTENS-1 model since more of its transmission lines were reinforced, thus allowing for more load to be served.  Based on the observation that line reinforcement is relatively inexpensive in our model, we have observed that the heuristic of reinforcing all lines performs relatively well. We also considered this approach as a comparison to the previous one.  The results of both Step (2) transmission approaches are given in Table \ref{tab:ph_results}.  Developing more sophisticated techniques for mapping stochastic CEP results from a low resolution model to a high resolution model is left for future work.  The results in Table \ref{tab:ph_results} demonstrate that most of the reduced models perform poorly on the stochastic CEP problem, with even KITTENS-1 resulting in 37\% error, corresponding to \$3.9B/year on this test system.  The most ``zonal'' model, KITTENS-100, resulted in 43\% error, corresponding to \$4.5B/year. We found that the higher costs of the reduced CEP models (as mapped to the baseline model) were largely due to higher capital costs (due to excess transmission reinforcement) but were also due to a mix of higher generation production costs and increased penalties for load shedding and RPS violation.

\begin{figure*}[t]
    \centering
    \includegraphics[width=\linewidth]{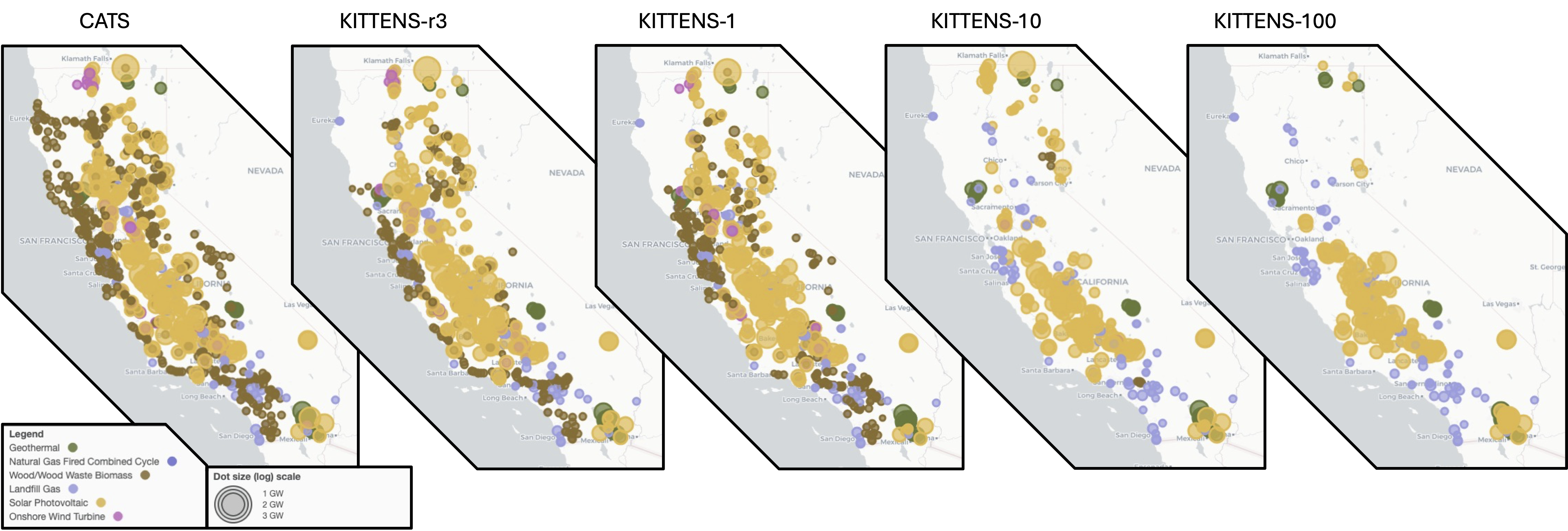}
    \caption{Comparison of different generation investment plans (as mapped to the CATS model) obtained from solving stochastic CEP problems of different spatial resolutions.}
    \label{fig:investment_plans}
\end{figure*}

Similar to the average across the single-day deterministic scenarios, the stochastic models at lower resolutions result in under-investment in new capacity.  For these simulations, the KITTENS-1 model results in the least amount of under-investment at 6\% and the KITTENS-100 model results in the most at 36\%.  Just as before, with the coarser models, we see large under-investment in new solar.  We also see under-investment in wind, wood biomass, natural gas, and storage in most cases.  The comparison in aggregate generation and storage investment decisions between the different models is given in Table \ref{tab:buildout_comparison_stochastic}.  We also visualize the different generation investment plans obtained from the different models (as mapped back to the CATS model) in Figure \ref{fig:investment_plans}. The more reduced models result in fewer sites for new generation and do not capture the value of many project sites in the full-resolution model, e.g., the KITTENS-100 model results in very little solar PV investment in northern California compared to the higher-resolution models.

The scope of this study was to look at CEP model coarsening for only \textit{one} test system and weather year, based on \textit{one} spatial aggregation heuristic.  Therefore, the results should not be taken as a generalization that is true for any power system, input data, and/or aggregation technique.  Even so, our results align with some existing work (on different test systems) that found solar under-investment in CEP models with lower spatial resolutions \cite{Jacobson_2024, krishnan_evaluating_2016}.  Unlike these papers, we also found under-investment in wind at lower spatial resolutions.  As noted in \cite{Jacobson_2024}, one cannot predict \textit{a priori} if a model will over- or under-invest in different types of generation at different spatial resolutions, and thus it is hard to generalize between models. Nonetheless, it is intuitive that lower spatial resolution CEP models distort investment decisions that rely on high spatial fidelity in their modeling (such as wind and solar) and also result in overall under-investment due to their ignorance of deliverability constraints.

\begin{table*}[t]
    \centering
        \caption{Solving 360-day Stochastic CEP at Varying Spatial Resolutions}
    \begin{tabular}{|l|l|l|l|l|l|l|l|l|}
    \hline
     &  \multicolumn{3}{c|}{\textit{\textbf{Step (1) via PH}}} & \multicolumn{5}{c|}{\textit{\textbf{Step (2) Results (Map A for Generation, Storage)}}}\\
    \hline
      \textbf{Case name}   & \textbf{PH solve} & \textbf{\# iters.} &\textbf{PH opt.} & \textbf{Transmission} &\textbf{Expected unserved} & \textbf{Loss of load} & \textbf{Achieved} &\textbf{Error}\\
         & \textbf{time (s)}&\textbf{of PH} &\textbf{gap}  & \textbf{approach}& \textbf{energy (EUE) } & \textbf{hours (LOLH)}  & \textbf{RPS} &\textbf{defined by (\ref{eqn:error_stochastic})}\\
         \hline\hline
        CATS & 72000*  & 57 & 5.16\% & N/A& 3.1 MWh & 2h / 8640h & 63.0\% & N/A\\ 
        \hline\hline
        KITTENS-r3 tightened & 72000* &89& 2.27\%  & Map solution&10.0 MWh & 9h / 8640h &61.2\% & \textbf{1.2\%}\\ 
& &&   & Reinforce all& 9.8 MWh & 8h / 8640h & 61.3\% & 36.3\% \\ 
\hline
        KITTENS-1 tightened & 72000*& 151 & 3.15\% &Map solution& 598.7 MWh & 1243h / 8640h & 59.2\% & 68.9\%\\ 
               &&&&Reinforce all & 9.9 MWh & 5h / 8640h & 61.6\% & 36.8\%\\ 
        \hline
        KITTENS-10 tightened &  2720 &41& 0.84\% & Map solution&513.6 MWh & 922h / 8640h & 55.2\% & 72.0\% \\ 
             &&&& Reinforce all & 45.6 MWh & 12h / 8640h & 58.8\% & 39.2\% \\ 
\hline
          KITTENS-100 tightened & 522 &21&  0.74\% &Map solution& 120.7 MWh &45h / 8640h & 54.4\% & 45.2\%\\ 
       &&&& Reinforce all & 52.7 MWh &15h / 8640h & 56.7\% & 43.0\%\\ 
         \hline
     \multicolumn{9}{r}{\textit{* = hit time limit}}\\
    \end{tabular}
    \label{tab:ph_results}
\end{table*}

\section{Conclusions}
This work presents a computational approach for verifying the assumption that zonal models can be used for CEP and for quantifying the impact of model coarsening.   While this work only presents one possible approach (of many) to CEP model coarsening, we are the first to demonstrate that rigorous computational validation on a high-fidelity nodal model can be used to either justify or reject a zonal approximation method. For our study, we introduce a heuristic approximation method for generating spatially aggregated models and consider several strategies for mapping results from these models back to a high-resolution nodal model.  
Using a synthetic-but-realistic California test case, we found that the proposed CEP models with low spatial resolution resulted in high errors and large amounts of generation under-investment (up to 41\%) when compared to the full-resolution model.  We found several CEP models with more limited spatial aggregation that were able to approximate the full-resolution model within reasonable tolerances. This work motivates ongoing research in both improved CEP-specific model reduction techniques as well as algorithmic advances to improve the tractability of nodal, high-fidelity CEP models.

\bibliographystyle{IEEEtran}
\bibliography{refs.bib}

\end{document}